\documentclass[12pt]{amsart}
\usepackage{amsmath,amssymb,amsbsy,amsfonts,amsthm,latexsym,
                        amsopn,amstext,amsxtra,euscript,amscd,mathrsfs,color,bm}

\usepackage{float}
\usepackage[english]{babel}
\usepackage{url}
\usepackage[colorlinks,linkcolor=blue,anchorcolor=blue,citecolor=blue,backref=page]{hyperref}

\usepackage{mathtools}
\usepackage{todonotes}
\usepackage[norefs,nocites]{refcheck}

\renewcommand*{\backref}[1]{}
\renewcommand*{\backrefalt}[4]{%
    \ifcase #1 (Not cited.)%
    \or        (p.\,#2)%
    \else      (pp.\,#2)%
    \fi} 

\newtheorem{theorem}{Theorem}
\newtheorem{lemma}[theorem]{Lemma}

\theoremstyle{definition}

\numberwithin{equation}{section}
\numberwithin{theorem}{section}
\numberwithin{table}{section}
\numberwithin{figure}{section}

\allowdisplaybreaks

\def\C{\mathbb{C}}
\def\G{\mathbb{G}}
\def\Z{\mathbb{Z}}

\def\fp{\mathfrak{p}}

 \def\mand{\qquad\mbox{and}\qquad}

\def\ff{\mathbf{f}}
\def\ov#1{{\overline{#1}}}

\def\fp{\mathfrak p}

\newcommand{\rA}{\ensuremath{\mathscr{A}}}

\def\({\left(}
\def\){\right)}

\newcommand{\g}{\gamma}

\newcommand{\Res}{{\mathrm{Res}}}
\newcommand{\h}{\mathrm{h}}
\newcommand{\di}{\mathrm{div}}
\newcommand\Gm{\G_{\mathrm{m}}}

\begin{document}

\title{Multiplicative dependence modulo subsets}

\author{Min Sha}
\address{School of Mathematical Sciences, South China Normal University, Guangzhou 510631, China}
\email{min.sha@m.scnu.edu.cn}

\subjclass[2010]{11N25, 11R04}

\keywords{Multiplicative dependence, multiplicative dependence modulo a subset}

\begin{abstract}
In this paper, 
we show that if the non-constant rational functions $f_1, \ldots, f_n\in K(x)$ over a number field $K$ 
   cannot multiplicatively generate a power of a linear fractional function, then
there are only finitely many elements $\alpha \in K$ such that $f_1(\alpha),\ldots,f_n(\alpha)$ are multiplicatively dependent modulo some subset `close' (with respect to the Weil height) to the division group of a finitely generated multiplicative subgroup of $K$. 
This improves some previous results. 
\end{abstract}

\maketitle


\section{Introduction}

\subsection{Motivation}
Let $\C$ be the field of complex numbers, and let $\Z$ be the ring of rational integers. 
Given  non-zero complex numbers $\alpha_1, \ldots, \alpha_n \in \C^*$, we say that they are \textit{multiplicatively dependent} if there exist integers $k_1,\ldots,k_n \in \Z$, not all zero, such that
\begin{equation*}
\alpha_1^{k_1}\cdots \alpha_n^{k_n} = 1;
\end{equation*}
and we say that they are \textit{multiplicatively dependent modulo $G$}, where $G$ is a subset of $\C^*$,
if there exist integers $k_1,\ldots,k_n \in \Z$, not all zero, such that
\begin{equation*}
\alpha_1^{k_1}\cdots \alpha_n^{k_n} \in G.
\end{equation*}
In addition, for an integer $r$ with $0 \le r \le n$, we say that $\alpha_1, \ldots, \alpha_n$ are 
multiplicatively dependent modulo $G$ of \textit{rank} $r$ if $r$ is the largest integer such 
that any $r$ of them are not multiplicatively dependent modulo $G$. 

Multiplicative dependence of algebraic numbers has been studied extensively 
from various aspects and still very actively;
see, for instance, \cite{AIS, BCMOS, BaS, BPV, BBGMOS, BBHOS, BHOS, BOSS, BMZ,  DS2016, DS18, KSSS, LvdP, Mello, OSSZ1, OSSZ2, PSSS, vdPL,  Stew, VZ, Young} and the references therein.
In \cite{BBGMOS}, the authors established a series of finiteness results about multiplicative dependence among values of rational functions modulo some subset which is somehow close to the division group of a finitely generated group $\Gamma$, denoted by $\Gamma^{\di}_\varepsilon$ (which is defined in the next section). 
In this paper, we want to continue this direction of study and improve some results in \cite{BBGMOS}.

The motivation also partly comes from the study of points on subvarieties of tori.
Let  $\Gm$ be  the multiplicative algebraic group over $\C$, that is $\Gm = \C^*$ endowed with the multiplicative group law. Let $H_d$ be the union of algebraic subgroups of $\Gm^n$ ($n \ge 2$) and with dimension at most $d$. 
Let $C$ be an absolutely irreducible curve in $\Gm^n$ defined over a number field 
and not contained in a translate of a proper algebraic subgroup of $\Gm^n$. 
In 1999, Bombieri, Masser and Zannier \cite{BMZ} proved that the points in $C \cap H_{n-1}$ are of bounded Weil height 
and the intersection $C \cap H_{n-2}$ is a finite set. 
Later, Maurin \cite{Mau} obtained the finiteness of $C \cap H_{n-2}$ under a weaker condition that 
$C$ is not contained in a proper algebraic subgroup of $\Gm^n$. 
Moreover, Maurin \cite{Mau11} showed that the intersection $C \cap \Gamma H_{n-2}$ is a finite set, 
where $\Gamma$ is a subgroup of $\Gm^n$ with finite rank. 
Recall that a proper algebraic subgroup of $\Gm^n$ is 
defined by finitely many equations of the form $x_1^{k_1} \cdots x_n^{k_n} = 1$ 
with unknowns $x_1, \ldots, x_n$ and integers $k_1, \ldots, k_n$ not all zero. 
If we let $G$ be the subgroup generated by all the coordinates 
of the points in $\Gamma$, then the coordinates of any point in $C \cap \Gamma H_{n-2}$ 
are multiplicatively dependent modulo $G$. 
Please see \cite{BMZ05, BMZ07, Hab, Mau11, Rem, Za} for some more works about intersections of subvarieties and algebraic subgroups in $\Gm^n$.



\subsection{Notations}
\label{subsection:notation}
Throughout the paper, we use the following notations:
\begin{itemize}
 \item[\textbullet]  $K$ is a number field.
 \item[\textbullet]  $\ov{K}$ is an algebraic closure of $K$.
 \item[\textbullet] $S$ is a finite set of places of $K$ containing all the infinite places.
 \item[\textbullet] $O_S$ is the ring of $S$-integers of $K$.
 \item[\textbullet] $R^*$ is the unit group of a ring $R$.
 \item[\textbullet]  $\Gamma$ is  a finitely generated subgroup of $K^*$.
 \item[\textbullet]  
The division group of $\Gamma$ is 
 $$\Gamma^{\di}:=\{\alpha \in \ov{K} : \alpha^m \in \Gamma \text{ for some integer } m \ge 1\}.$$
  \item[\textbullet] For $\varepsilon > 0$, $$\Gamma^{\di}_\varepsilon := \{\alpha \beta : \alpha \in \Gamma^{\di}, \beta \in \ov{K}^* \text{ with } \h(\beta) \leq \varepsilon \},$$
where $\h(\cdot)$ stands for the absolute logarithmic Weil height function.
 \end{itemize}


In addition, let $M_K$ be the set of places of $K$, $M_K^\infty$ the set of infinite places of $K$,
and $M_K^0 = M_K \setminus M_K^\infty$.

\subsection{Main results}

Recall that a linear fractional function over $\C$ is of the form $(ax + b)/(cx+d)$ with variable $x$,
 constants $a, b, c, d \in \C$ and $cx+d \ne 0$. 

We say that the rational functions $f_1,\ldots,f_n\in\C(x)$ \textit{multiplicatively generate} a rational function $g$
if there exist integers $k_1,\ldots,k_n \in \Z$, not all zero, such that
$$
f_1^{k_1} \cdots f_n^{k_n} = g.
$$

\begin{theorem}
\label{thm:multdep}
Let $f_1, \ldots, f_n \in K(x)$ be non-constant rational functions.
Assume that they cannot multiplicatively generate a power of a linear fractional function. 
Then, for every $\varepsilon > 0$ there are only finitely many elements $\alpha \in K$ such that $f_1(\alpha), \ldots, f_n(\alpha)$ are multiplicatively dependent modulo $\Gamma^{\di}_\varepsilon$.
\end{theorem}

We remark that 
some finiteness results similar as the above theorem and under more restrictive conditions were given in 
\cite[Theorems 1.1 and 1.2]{BBGMOS}.
Especially, the case when $n=2$ and both $f_1$ and $f_2$ are polynomials was established 
in \cite[Theorem 1.5]{BBGMOS}.

We also remark that in Theorem~\ref{thm:multdep} the condition ``they cannot multiplicatively generate a power of a linear fractional function" can not be removed in general. For example, let $n=1$ and $f_1(x) = a(x-b)/(x-c)$ with $a, b, c \in K$ and $a(b-c)\ne 0$; and then choosing any $\beta \in \Gamma^{\di}_\varepsilon \cap K$ with $\beta \ne a$ and puting $\alpha = (ab-c\beta)/(a-\beta) \in K$, one can check directly that $f_1(\alpha) = \beta \in \Gamma^{\di}_\varepsilon$.

Next, we somehow relax the condition on the rational functions $f_1, \ldots, f_n$ but enhance the condition on multiplicative dependence. 

For any rational function $f \in K(x)$, the numerator and denominator of $f$ are meant to be two polynomials $g, h \in K[x]$, 
respectively, such that $f= g / h$, $\gcd(g,h)=1$ and $h$ is monic.   
Moreover, we say that $u \in K[x]$ is a factor of $f$ if $u$ divides the numerator or the denominator of $f$.

\begin{theorem}  
\label{thm:multdep2}
Let $f_1, \ldots, f_n  \in K(x)$ be non-constant rational functions $(n \ge 2)$. 
Assume that at least one of them has one irreducible factor over $K$ which has degree greater than 1 
and is not a factor of the other $n-1$ rational functions.  
Then, for every $\varepsilon > 0$, the following set is finite: 
\begin{align*}
\Big\{\alpha \in K: \ & f_1(\alpha),  ..., f_n(\alpha) \textrm{ are multiplicatively} \\
& \qquad\qquad \textrm{dependent mod $\Gamma^{\di}_\varepsilon$ of rank $n-1$}\Big\}.
\end{align*}
\end{theorem}



\section{Preliminaries}
In this section, we collect some results used later on. 

We define the set $\rA(K,H)$ as the set of nonzero elements in the algebraic number
field $K$ of height at most $H$, that is,
$$
\rA(K,H) =\left\{\alpha \in  K^*:~ \h(\alpha)\le H \right\}.
$$
We note that by Northcott's Theorem the set $ \rA(K, H)$  is a finite set.

We need the following result from \cite[Theorem 2.1]{OstShp}.

\begin{lemma}
\label{lem:K H}
Let $\{g_1, \ldots, g_r\}$ be a set of generators of $\Gamma$, which minimises $H = \max_{i=1, \ldots, r} \h( g_i)$. Then, for every
$\varepsilon> 0$, we have
$$
K^*\cap \Gamma^{\di}_\varepsilon \subseteq \left \{\beta \eta :~ 
\beta \in \Gamma, \eta \in \rA(K, \varepsilon  + rH)\right\}.
$$
\end{lemma}


As usual, for any non-constant rational function $f \in K(X)$, the degree of $f$ is defined to be the maximum of the degrees of its numerator and denominator.  
The following result is \cite[Theorem 1.2 (a)]{BOSS}.

\begin{lemma}
\label{lem:KGa}
Let $f \in K(x)$ be a rational function of degree $d \ge 2$.
Assume that $f$ is not of the form $a(x-b)^d$ or $a(x-b)^d / (x-c)^d$ with $a,b,c \in K$, $a(b-c)\ne0$, and $d \in \Z$.
Then, the set $\{\alpha \in K: \, f(\alpha) \in \Gamma \}$ is finite.
\end{lemma}


\section{Proofs of the main results}

\subsection{Preliminary discussion}
\label{sec:pre}

Let~$S_\Gamma$ be the following set of places of $K$:
\begin{equation*}
  S_\Gamma := M_K^\infty \cup \big\{v\in M_K^0 : \,  \text{$v(\gamma) \ne 0$ for some $\g\in\Gamma$} \big\},
\end{equation*}
where, as usual, $v(\gamma)$ means the additive valuation of $v$ at $\gamma$.
Note that the set~$S_\Gamma$ is finite, since~$\Gamma$ is finitely generated.

As usual, we say that a polynomial
$$
f(x) =  a_d x^d+ \cdots + a_{1}x+ a_0 \in K[x]
$$
has bad reduction at~$v\in M_{K}^0$
if either $v(a_i)<0$ for some~$i=0,1,\ldots,d-1$, or $v(a_d)\ne0$; otherwise
we say it has good reduction at $v$. 

Moreover, for a rational function $f(x) \in K(x)$, we say that $f$ has bad reduction at~$v\in M_{K}^0$
if either the numerator or the denomenator of $f$  has bad reduction at~$v$. 

For the rational functions $f_{1},   \ldots, f_{n} \in K(x)$ given in the theorem, 
let $g_1, \ldots, g_m$ be all the distinct monic irreducible factors (over $K$) in the numerators and denominators of  $f_{1},  \ldots, f_{n}$. 
Then, for each $f_i, 1\le i \le n$, we can write 
\begin{equation}  \label{eq:figj}
f_i = a_i \prod_{j=1}^{m} g_j^{e_{ij}},\quad a_i\in K^*,
\end{equation}
for some integers $e_{i1}, \ldots, e_{im}$ (which are not all zero). 
Besides, we denote 
$$
\ff=(f_1,\ldots,f_n). 
$$

Now, we define the finite set $S_{\ff,\Gamma}$ to be the union of the following three parts of places: 
the first part is $S_\Gamma$, the second part is the subset of places $v\in M_K^0$ such that 
at least one of $f_1, \ldots, f_n, g_1, \ldots, g_m$ has bad reduction at $v$, 
and the third part is the subset of places $v\in M_K^0$ such that 
at least one of  the resultants $\Res(g_i, g_j)$, $1 \le i \ne j \le m$, has non-zero valuation at $v$. 


Moreover, for any $\varepsilon >0$, by Lemma~\ref{lem:K H} we have that every element in 
the intersection $K^*\cap \Gamma^{\di}_\varepsilon$ is of the form $\beta\eta$
, where $\beta \in \Gamma$ and $\eta \in K^*$ with $\h(\eta)\ll_{\varepsilon,\Gamma} 1$. 
Since $\eta \in K^*$ is of bounded height depending only on $\varepsilon$ and $\Gamma$, by Northcott's theorem there are only finitely many such $\eta$. Thus we can enlarge the set $S_{\ff,\Gamma}$ to include all prime ideals that divide these finitely many elements $\eta$. 
We denote the new set by $S_{\ff,\Gamma,\varepsilon}$, and we also note that $S_{\ff,\Gamma,\varepsilon}$ is still a finite set.

By the construction of the set $S_{\ff,\Gamma, \varepsilon}$,  we have
$$
\Gamma \subseteq  O_{S_{\ff,\Gamma, \varepsilon}}^* \mand  K^*\cap \Gamma^{\di}_\varepsilon \subseteq  O_{S_{\ff,\Gamma, \varepsilon}}^*.
$$
Hence, it suffices to prove the main results by replacing $ \Gamma^{\di}_\varepsilon$ with $O_{S_{\ff,\Gamma,\varepsilon}}^*$.


\subsection{Proof of Theorem~\ref{thm:multdep}}

Following the above discussions, we only need to prove the theorem by replacing $ \Gamma^{\di}_\varepsilon$ with $O_{S_{\ff,\Gamma,\varepsilon}}^*$.

 If $n = 1$, since $f_1$ is not a power of any linear fractional function and $O_{S_{\ff,\Gamma,\varepsilon}}^*$ is a finitely generated subgroup, we see that applying  Lemma~\ref{lem:KGa} to $f_1$ and $O_{S_{\ff,\Gamma,\varepsilon}}^*$ 
 gives the desired finiteness result.

In the sequel, we assume that $n \ge 2$, and that the result is valid for $n-1$, in order to apply an induction.

Let $\alpha \in K$ be such that $ f_1(\alpha), \ldots, f_n(\alpha) \not\in O_{S_{\ff,\Gamma,\varepsilon}}^*$ and there exist integers $k_1,..., k_n$, not all zero such that
$
f_1(\alpha)^{k_1} \cdots  f_n(\alpha)^{k_n} \in O_{S_{\ff,\Gamma,\varepsilon}}^*.
$
We write
\begin{equation}
\label{eq:multrs2}
f_1(\alpha)^{k_1} \cdots  f_n(\alpha)^{k_n}= \gamma,\qquad  \gamma  \in O_{S_{\ff,\Gamma,\varepsilon}}^*.
\end{equation}
By induction, we can assume $k_1 \cdots k_n \ne 0$.


Fix a prime ideal $\fp \not\in S_{\ff,\Gamma,\varepsilon}$ in $ K$. 
If $v_\fp(\alpha)<0$, then using the ultrametric inequality of non-Archimedean valuations 
we have 
\begin{equation*}
\label{eq:ord gj}
   v_\fp(g_{j}(\alpha)) = d_j v_\fp(\alpha), \quad d_j = \deg(g_j), \quad\text{for  $j=1,\ldots,m$}.
\end{equation*}
Then, in view of \eqref{eq:figj}, for each $i=1, \ldots, n$ we get 
\begin{equation}   \label{eq:vpfi}
 v_\fp(f_i(\alpha)) = \sum_{j=1}^{m} e_{ij} v_{\fp}(g_j(\alpha)) = \sum_{j=1}^{m} e_{ij} d_j v_{\fp}(\alpha), 
\end{equation}
where we also use the fact $v_\fp(a_i)=0$ due to the constructions of $S_{\ff,\Gamma}$ and $ S_{\ff,\Gamma,\varepsilon}$.
Combining this with \eqref{eq:multrs2}, we obtain (noticing $v_\fp(\gamma)=0$ due to $\gamma\in O_{S_{\ff,\Gamma,\varepsilon}}^*$)
\begin{equation}
\label{eq:exponent1}
k_1 \sum_{j=1}^{m} d_j e_{1j}  + \cdots + k_n \sum_{j=1}^{m} d_j e_{nj} = 0.
\end{equation}

Otherwise, if $v_\fp(\alpha) \ge 0$,  then due to the constructions of $S_{\ff,\Gamma}$ and $ S_{\ff,\Gamma,\varepsilon}$, 
we have $v_\fp(g_j(\alpha)) \ge 0$ for each $j=1, \ldots, m$, and 
moreover, at most one of them is positive. 
Indeed, if  $v_\fp(g_1(\alpha)) > 0$ and  $v_\fp(g_2(\alpha)) > 0$ without loss of generality, 
then this implies that $v_\fp(\Res(g_1,g_2))>0$ (notice that, since $g_1,g_2 \in O_{S_{\ff,\Gamma,\varepsilon}}[X]$ and $g_1$ and $g_2$ do not have common roots, we have $\Res(g_1,g_2)\in O_{S_{\ff,\Gamma,\varepsilon}}$ and $\Res(g_1,g_2) \ne 0$).
By our construction of the set $S_{\ff,\Gamma,\varepsilon}$, this implies that $\fp \in S_{\ff,\Gamma,\varepsilon}$,
which is a contradiction with the choice of $\fp$ above. 

Hence, if $v_\fp(\alpha) \ge 0$ and $v_\fp(g_j(\alpha))>0$ for some $j$, 
then for any $t \in \{1, \ldots, m\}$ with $t \ne j$  we 
have $v_\fp(g_t(\alpha))=0$; and so, considering valuations in \eqref{eq:multrs2} and using \eqref{eq:figj}, we obtain 
\begin{equation}
\label{eq:exponent2}
k_1  e_{1j}  + \cdots + k_n e_{nj} = 0.
\end{equation}

Now, without loss of generality,
by varying the prime ideal  $\fp \not\in S_{\ff,\Gamma,\varepsilon}$ in $ K$ 
 we can assume that 
the identity \eqref{eq:exponent1} and for each $j=1, \ldots, m$ the identity \eqref{eq:exponent2}  can appear. 
Then, $(k_1, \ldots, k_n)$ is in fact an integer solution of the system of linear Diophantine equations 
in \eqref{eq:exponent1} and \eqref{eq:exponent2} (totally having $m+1$ equations).  
Since the integer solutions of this system form a submodule of 
the free $\Z$-module $\Z^n$  and $\Z$ is a principal ideal domain, 
we know that they in fact form a free $\Z$-module. 

Hence, there exists a basis of the integer solutions  of this system, say,
$$
(b_{l1}, \ldots, b_{ln}), \quad l = 1, \ldots, r. 
$$
Here we emphasize that this basis is independent of $\alpha$.
Then, the integers $k_1, \ldots, k_n$ can be expressed as
$$
k_i = \sum_{l=1}^{r} c_l b_{li},  \quad i=1,\ldots,n,
$$
for some integers $c_1, \ldots, c_{r}$.

Now, we let 
\begin{equation*}  \label{eq:F}
F(x) = \prod_{i=1}^{n} f_i(x)^{b_{1i}}, 
\end{equation*}
where the exponent vector $(b_{11}, \ldots, b_{1n})$ is non-zero by its choice above 
and is also an integer solution of that system. 

For any $\fp \not\in S_{\ff,\Gamma,\varepsilon}$, 
if $v_\fp(\alpha)<0$, then \eqref{eq:vpfi} holds and  
we obtain 
\begin{align*}
v_\fp(F(\alpha))= \sum_{i=1}^{n}b_{1i} v_\fp(f_i(\alpha)) = 
v_\fp(\alpha) \sum_{i=1}^{n}b_{1i} \sum_{j=1}^{m} d_j e_{ij} = 0, 
\end{align*}
where the last equality follows from the fact that 
$(b_{11}, \ldots, b_{1n})$  is a solution to \eqref{eq:exponent1}. 
Otherwise, if $v_\fp(\alpha) \ge 0$, then as before we have that 
the valuations of $g_1(\alpha), \ldots, g_m(\alpha)$ at $\fp$ are all non-negative and 
 at most one of them is positive; say $v_{\fp}(g_j(\alpha))>0$ for some $j$, and then we get 
$$
v_\fp(F(\alpha))= \sum_{i=1}^{n}b_{1i} v_\fp(f_i(\alpha)) = v_\fp(g_j(\alpha)) \sum_{i=1}^{n}b_{1i} e_{ij} = 0,  
$$
since $(b_{11}, \ldots, b_{1n})$  is also a solution to \eqref{eq:exponent2}. 
Hence, we conclude that for any prime ideal $\fp \not\in S_{\ff,\Gamma,\varepsilon}$, 
we have $v_\fp(F(\alpha))=0$. 
So, we have $F(\alpha) \in O_{S_{\ff,\Gamma,\varepsilon}}^*$.
Now, the desired result follows directly from Lemma~\ref{lem:KGa} 
(which we can apply, since by assumption $F$ is not a power of a linear fractional function). 
This completes the proof.

\subsection{Proof of Theorem~\ref{thm:multdep2}}

Let $\alpha \in K$ is in the set described in the theorem. 
Then, choosing the set of places $S_{\ff,\Gamma,\varepsilon}$ as before, 
we know that there exist non-zero integers $k_1,..., k_n$ such that
$$
f_1(\alpha)^{k_1} \cdots  f_n(\alpha)^{k_n}= \gamma 
$$
for some element $\gamma \in O_{S_{\ff,\Gamma,\varepsilon}}^*$. 

Applying the same arguments as the above, we deduce that there exists a basis of the integer solutions  
of a system of linear Diophantine equations, say
$$
(b_{l1}, \ldots, b_{ln}), \quad l = 1, \ldots, r, 
$$
such that the integers $k_1, \ldots, k_n$ can be expressed as
$$
k_i = \sum_{l=1}^{r} c_l b_{li},  \quad i=1,\ldots,n,
$$
for some integers $c_1, \ldots, c_{r}$. 
We emphasize that this basis is independent of $\alpha$.

By assumption and without loss of generality, we can assume that 
$f_1$ has one irreducible factor over $K$ which has degree greater than 1 
and is not a factor of $f_2 \cdots f_n$.  
Since $k_1 \ne 0$, we know that there exists an integer $l \in \{1, \ldots, r\}$ such that $b_{l1} \ne 0$. 
Then, we define 
\begin{equation*}  
G(x) = \prod_{i=1}^{n} f_i(x)^{b_{li}}, 
\end{equation*}
where $b_{l1} \ne 0$ and the exponent vector $(b_{l1}, \ldots, b_{ln})$  is an integer solution of that system. 
Then, as before we have $G(\alpha) \in O_{S_{\ff,\Gamma,\varepsilon}}^*$.
Now, the desired result follows directly from Lemma~\ref{lem:KGa} 
(which we can apply, since by assumption $G$ is not a power of a linear fractional function). 
This completes the proof.

\section*{Acknowledgement}

The author is grateful to Alina Ostafe for her valuable comments on an earlier version of this paper. 
The research was supported by the Guangdong Basic and Applied Basic Research Foundation (No. 2025A1515010635).

\end{document}